\documentclass[a4,12pt]{elsart}
\usepackage{ifpdf}
\usepackage{amsmath,amsfonts,latexsym,amssymb}
\usepackage{graphicx,indentfirst}
\ifpdf

\usepackage[%
  pdftitle={Instructions for use of the document class
    elsart},%
  pdfauthor={Simon Pepping},%
  pdfsubject={The preprint document class elsart},%
  pdfkeywords={instructions for use, elsart, document class},%
  pdfstartview=FitH,%
  bookmarks=true,%
  bookmarksopen=true,%
  breaklinks=true,%
  colorlinks=true,%
  linkcolor=blue,anchorcolor=blue,%
  citecolor=blue,filecolor=blue,%
  menucolor=blue,pagecolor=blue,%
  urlcolor=blue]{hyperref}
\else
\usepackage[%
  breaklinks=true,%
  colorlinks=true,%
  linkcolor=blue,anchorcolor=blue,%
  citecolor=blue,filecolor=blue,%
  menucolor=blue,pagecolor=blue,%
  urlcolor=blue]{hyperref}
\fi

\makeatletter
\def\elsartstyle{%
    \def\normalsize{\@setfontsize\normalsize\@xiipt{14.5}}
    \def\small{\@setfontsize\small\@xipt{13.6}}
    \let\footnotesize=\small
    \def\large{\@setfontsize\large\@xivpt{18}}
    \def\Large{\@setfontsize\Large\@xviipt{22}}
    \skip\@mpfootins = 18\p@ \@plus 2\p@
    \normalsize
               }
\@ifundefined{square}{}{} \makeatother

\newtheorem{theorem}{Theorem}
\newtheorem{lemma}[theorem]{Lemma}
\newtheorem{corollary}[theorem]{Corollary}
\newtheorem{proposition}[theorem]{Proposition}
\newcommand{\fd}{\mathbb{F}}
\newcommand{\z}{\mathbb{Z}}

\newcommand\proof{{\sc Proof.} \enspace}

\begin{document}
\begin{minipage}{\textwidth}
\elsartstyle
\parskip 12pt
\renewcommand{\thempfootnote}{\fnsymbol{mpfootnote}}
\leftskip=2pc
\begin{center}
{\LARGE Infinite Families of Recursive Formulas\\
\ Generating Power Moments of Kloosterman Sums: $O^{+}(2n,2^{r})$
Case\par}
\large
Dae San Kim\\[12pt]
\small\itshape Department of Mathematics, Sogang University, Seoul
121-742, Korea
\end{center}

\bigskip
\leftskip=0pt

\hrule\vskip 8pt
\begin{small}
{\bfseries Abstract}
\parindent 1em

In this paper, we construct four infinite families of binary linear
codes  associated with double cosets with respect to certain maximal
parabolic subgroup of the orthogonal group $O^{+}(2n,2^{r})$. Here
$q$ is a power of two. Then we obtain two infinite families of
recursive formulas for the power moments of Kloosterman sums and
those of 2-dimensional Kloosterman sums in terms of the frequencies
of weights in the codes. This is done via Pless power moment
identity and by utilizing the explicit expressions of exponential
sums over those double cosets  related to the evaluations of ``Gauss
sums" for the orthogonal  groups $O^{+}(2n,2^{r})$.

\end{small}
\vskip 10pt\hrule

\leftskip=0pt

\vspace{24pt}
\renewcommand{\thempfootnote}{\astsymbol{mpfootnote}}
\footnotetext[1]{Corresponding author.} \setbox0=\hbox{\footnotesize
1} \edef\thempfootnote{\hskip\wd0} \footnotetext[0]{\textit{Email
adress:} dskim@sogang.ac.kr
  (Dae San Kim).}

\section {Introduction}

Let  $\psi$ be a nontrivial additive character of the finite field
$\fd_q$ with $q = p^r$ elements ( $p$ a prime), and let $m$ be a
positive integer. Then the $m$-dimensional  Kloosterman sum
$K_m(\psi;a)$(\cite{RH}) is defined by

\begin{equation*}
 K_{m}(\psi;a)=\sum_{\alpha_{1},\cdots,\alpha_{m} \in
\fd_{q}^{*}}\psi(\alpha_{1}+\cdots+\alpha_{m}+a\alpha_{1}^{-1}\cdots\alpha_{m}^{-1})\\
(a \in \fd_{q}^{*}).
\end{equation*}\\

\end{minipage}
\bigskip

In particular, if $m=1$, then $K_{1}(\psi;a)$ is simply denoted by
$K(\psi;a)$, and is called the Kloosterman sum. The Kloosterman sum
was introduced in 1926 to give an estimate for the Fourier
coefficients of modular forms (cf. \cite{HDK}, \cite{JH}). It has
also been studied to solve  various problems in coding theory and
cryptography over finite fields of characteristic two (cf.
\cite{PTV}, \cite{HPT}).

For each nonnegative integer $h$, by $MK_{m}(\psi)^{h}$ we will
denote the $h$-th moment of the $m$-dimensional Kloosterman sum
$K_{m}(\psi;a)$. Namely, it is given by
\begin{equation*}
 MK_{m}(\psi)^{h}=\sum_{a \in \mathbb{F}_{q}^{*}}K_{m}(\psi;a)^{h}.
 \end{equation*}

If $\psi=\lambda$ is the canonical additive character of
$\mathbb{F}_{q}$, then $MK_{m}(\lambda)^{h}$ will be simply denoted
by $MK_{m}^{h}$. If further $m=1$, for brevity $MK_{1}^{h}$ will be
indicated by $MK^{h}$.

Explicit computations on power moments of Kloosterman sums were
begun with the paper \cite{HS} of Sali\'{e} in 1931, where he
showed, for any odd prime $q$,
\begin{equation*}
MK^h=q^2M_{h-1}-(q-1)^{h-1}+2(-1)^{h-1}(h \geq 1).
\end{equation*}

Here $M_{0}=0$, and for $h \in \mathbb{Z}_{>0}$,

\begin{equation*}
M_{h}=|\{(\alpha_{1}, \cdots,\alpha_{h})\in(\mathbb{F}_{q}^{*})^{h}
| \sum_{j=1}^{h} \alpha_{j}=1=\sum_{j=1}^{h}\alpha_j^{-1}\}|.
\end{equation*}

For $q=p$ odd prime, Sali\'{e} obtained  $MK^{1}$, $MK^{2}$, $MK^{3}
$, $MK^{4}$ in \cite{HS} by determining $M_{1},M_{2},M_{3}$.
$MK^{5}$ can be expressed in terms of the $p$-th eigenvalue for a
weight 3 newform on $\Gamma_{0}$(15) (cf. \cite{RL}, \cite{CJM}).
$MK^{6}$ can be expressed in terms of the $p$-th eigenvalue for a
weight 4 newform on $\Gamma_{0}$(6) (cf. \cite{KJBD}). Also, based
on numerical evidence, in \cite{RJE}. Also, based on numerical
evidence, in \cite{RJE} Evans was led to propose a conjecture which
expresses $MK^7$ in terms of Hecke eigenvalues for a weight 3
newform on $\Gamma_0 (525)$ with quartic nebentypus of conductor
105. For more details about this brief history of explicit
computations on power moments of Kloosterman sums, one is referred
to Section IV of \cite{D2}.

From now on, let us assume that $q=2^{r}$. Carlitz\cite{L1}
evaluated $MK^{h}$ for $h \leq 4$. Recently, Moisio was able to find
explicit expressions of $MK^{h}$, for the other values of $h$ with
$h \leq 10$ (cf.\cite{M1}). This was done, via Pless power moment
identity, by connecting moments of Kloosterman sums and the
frequencies of weights in the binary Zetterberg code of length
$q+1$, which were known by the work of Schoof and Vlugt in
\cite{RM}.

In \cite{D2}, the binary linear codes $C(SL(n,q))$ associated with
finite special linear groups $SL(n,q)$ were constructed when $n,q$
are both powers of two. Then obtained was a recursive formula for
the power moments of multi-dimensional Kloosterman sums in terms of
the frequencies of weights in $C(SL(n,q)$. In particular, when
$n=2$, this gives a recursive formula for the power moments of
Kloosterman sums. Also, in order to get recursive formulas for the
power moments of Kloosterman and 2-dimensional Kloosterman sums, we
constructed in \cite{D3} three binary linear codes $C(SO^{+}(2,q))$,
$C(O^{+}(2,q) )$, $C(SO^{+}(4,q))$, respectively associated with
$SO^{+}(2,q)$, $O^{+} (2,q)$, $SO^{+}(4,q)$, and in \cite{D4} three
binary linear codes $C(SO^{-}(2,q))$, $C(O^{-}(2,q))$,
$C(SO^{-}(4,q))$, respectively associated with $SO^{-}(2,q)$,
$O^{-}(2,q)$, $SO^{-}(4,q)$. All of these were done via Pless power
moment identity and by utilizing our previous results on explicit
expressions of Gauss sums for the stated finite classical groups.
Still, in all, we had only a handful of recursive formulas
generating power moments of Kloosterman and 2-dimesional Kloosterman
sums.

In this paper, we will be able to produce two infinite families of
recursive formulas generating power moments of Kloosterman sums and
two those of 2-dimensional Kloosterman sums. To do that, we
construct four infinite families of binary linear codes
$C(DC_{1}^{+}(n,q))$ $(n=2,4,\cdots)$, $C(DC_{1}^{-}(n,q))$
$(n=1,3,\cdots)$, both associated with $P^{+} \sigma_{n-1}^{+}
P^{+}$, and $C(DC_{2}^{+} (n,q))$ $(n=2,4, \cdots)$,
$C(DC_2^-(n,q))$ $(n=3,5,\cdots)$, both associated with $P^+
\sigma_{n-2}^+ P^+$, with respect to the maximal parabolic subgroup
$P^+=P^+(2n,q)$ of the orthogonal group $ O^+(2n,q)$, and express
those power moments in terms of the frequencies of weights in each
code. Then, thanks to our previous results on the explicit
expressions of exponential sums over those double cosets related to
the evaluations of ``Gauss sums" for the orthogonal  groups
$O^+(2n,q)$ [9], we can express the weight of each codeword in the
duals of the codes in terms of Kloosterman or 2-dimensional
Kloosterman sums. Then our formulas will follow immediately from the
Pless power moment identity. Similarly to these, in [7], we obtained
infinite families of recursive formulas for power moments of
Kloosterman and 2-dimensional Kloosterman sums by constructing
binary codes associated with double cosets with respect to certain
maximal parabolic subgroup of the symplectic group $Sp(2n,q)$

Theorem \ref{A} in the following(cf. (\ref{a9}), (\ref{a10}),
(\ref{a12})-(\ref{a14})) is the main result of this paper.
Henceforth, we agree that the binomial coefficient $\binom{b}{a}=0$,
if $a > b$ or $a<0$. To simplify notations, we introduce the
following ones which will be used throughout this paper at various
places.

\begin{equation}\label{a1}
A_{1}^{+}(n,q)=q^{\frac{1}{4}(5n^{2}-6n)}\left[ \substack{n \\ 1}
 \right]_q \prod_{j=1}^{n/2}(q^{2j-1}-1),
\end{equation}
\begin{equation}\label{a2}
B_{1}^{+}(n,q)=q^{\frac{1}{4}(n-2)^{2}}\prod_{j=1}^{n/2}(q ^{2j}-1),
\end{equation}
\begin{equation}\label{a3}
A_{2}^{+}(n,q)=q^{\frac{1}{4}(5n^{2}-6n)} \left[ \substack{n \\ 2}
 \right]_q \prod_{j=1}^{(n-2)/2}(q^{2j-1}-1),
\end{equation}
\begin{equation}\label{a4}
B_{2}^{+}(n,q)=q^{\frac{1}{4}(n^2-8n+12)}(q^{n-1}-1)(q^{n}-1)
\prod_{j=1}^{(n-2)/2}(q^{2j}-1),
\end{equation}
\begin{equation}\label{a5}
A_{1}^{-}(n,q)=q^{\frac{1}{4}(5n^2-4n-1)} \left[ \substack{n \\ 1}
\right]_q \prod_{j=1}^{(n-1)/2} (q^{2j-1}-1),
\end{equation}
\begin{equation}\label{a6}
B_{1}^{-}(n,q)=q^{\frac{1}{4}(n^{2}-6n + 5)}(q^n-1)
\prod_{j=1}^{(n-1)/2} (q^{2j}-1),
\end{equation}
\begin{equation}\label{a7}
A_{2}^{-}(n,q)=q^{\frac{1}{4}(5n^2-8n+3)}\left[ \substack{n \\
2}\right]_q \prod_{j=1}^{(n-1)/2}(q^{2j-1}-1),
\end{equation}
\begin{equation}\label{a8}
B_{2}^{-}(n,q)=q^{\frac{1}{4}(n-3)^2}(q^n-1)\prod_{j=1}^{(n-1)/2}(q^{2j}-1).
\end{equation}

From now on, it is assumed that either $+$ signs or $-$
 signs are
chosen everywhere, whenever $\pm$ signs appear. Henceforth we agree
that the binomial coefficient $\binom{b}{a}=0$, if $a>b$ and $a<0$.

\begin{theorem}:\label{A}
 Let $q=2^{r}$. Then, with the notations in (\ref{a1})-(\ref{a8}), we have the following.
next line With + signs everywhere for $\pm$ signs,  we have a
recursive formula generating power moments of Kloosterman sums over
$ \fd_{q}$, for each $n \geq 2$ even and all $q$. Also, with $-$
signs everywhere for $\pm$ signs, we have such a formula, for either
each $n \geq 3$ odd and all $q$, or $n=1$ and $q \geq 8$:

\begin{align}\label{a9}
  \begin{split}
 &MK^h=\sum_{l=0}^{h-1}(-1)^{h+l+1}\binom{h}{l}B_{1}^{\pm}(n,q)^{h-l}MK^{l}\\
      &+qA_{1}^{\pm}(n,q)^{-h}\sum_{j=0}^{min \{N_{1}^{\pm}(n,q),h\}}(-1)^{h+j}C_{1,j}^{\pm}(n,q)\sum_{t=j}^{h}t!S(h,t)2^{h-t}\binom{N_{1}^{\pm}(n,q)-j}{N_{1}^{\pm}(n,q)-t}\\
      & (h=1,2,\cdots),
  \end{split}
\end{align}

where
$N_{1}^{\pm}(n,q)=|DC_{1}^{\pm}(n,q)|=A_{1}^{\pm}(n,q)B_{1}^{\pm}(n,q)$,
and $\{C_{1,j}^{\pm}(n,q)\}_{j=0}^{N_{1}^{\pm}(n,q)}$ is the weight
distribution of $C(DC_{1}^{\pm}(n,q))$ given by

\begin{align}\label{a10}
  \begin{split}
 &C_{1,j}^{\pm}(n,q)=\sum_{}^{}\binom{q^{-1}A_{1}^{\pm}(n,q)(B_{1}^{\pm}(n,q)+1)}{\nu_{0}} \\
 &\times \prod_{tr(\beta^{-1})=0} \binom{q^{-1}A_{1}^{\pm}(n,q)(B_{1}^{\pm}(n,q)+q+1)}{\nu_{\beta}}\prod_{tr(\beta^{-1})=1}\binom{q^{-1}A_{1}^{\pm}(n,q)(B_{1}^{\pm}(n,q)-q+1)}{\nu_{\beta}}.
  \end{split}
\end{align}
Here the sum is over all the sets of nonnegative integers
$\{\nu_{\beta}\}_{\beta \in \mathbb{F}_{q}}$ satisfying $
\sum_{\beta \in \mathbb{F}_{q}} \nu_{\beta}=j$ and $\sum_{\beta \in
\mathbb{F}_{q}}^{} \nu_{\beta} \beta =0$. In addition, $S(h,t)$ is
the Stirling number of the second kind defined by
\begin{equation}\label{a11}
S(h,t)= \frac{1}{t!} \sum_{j=0}^{t}(-1)^{t-j}\binom{t}{j}j^{h}.
\end{equation}

(b) With $+$ signs everywhere for $\pm$ signs,  we have recursive
formulas generating power moments of 2-dimensional Kloosterman sums
over $\fd_{q}$ and even power moments of Kloosterman sums over
$\fd_{q}$, for each even $ n \geq 2$ and all $q \geq 4$. Also, with
$-$ signs everywhere for $\pm$ signs, we have such formulas, for
each $n \geq 3$ odd and $q \geq 4$:

\begin{align}\label{a12}
  \begin{split}
 &MK_2^h=\sum_{l=0}^{h-1}(-1)^{h+l+1}\binom{h}{l}(B_{2}^{\pm}(n,q)-q^2)^{h-l}MK_2^l \\
 &+qA_{2}^{\pm}(n,q)^{-h}\sum_{j=0}^{min \{ N_{2}^{\pm}(n,q),h \}}(-1)^{h+j}C_{2,j}^{\pm}(n,q)\sum_{t=j}^{h}t!S(h,t)2^{h-t} \binom{N_{2}^{\pm}(n,q)-j}{N_{2}^{\pm}(n,q)-t}\\
 &( h=1,2,\cdots),
  \end{split}
\end{align}
and
\begin{align}\label{a13}
  \begin{split}
 &MK^{2h}=\sum_{l=0}^{h-1}(-1)^{h+l+1}\binom{h}{l}(B_{2}^{\pm}(n,q)-q^2+q)^{h-l}MK^{2l}\\
 &+qA_{2}^{\pm}(n,q)^{-h} \sum_{j=0}^{min \{N_{2}^{\pm}(n,q),h \}}(-1)^{h+j}C_{2,j}^{\pm}(n,q)\sum_{t=j}^{h}t!S(h,t)2^{h-t}\binom{N_{2}^{\pm}(n,q)-j}{N_{2}^{\pm}(n,q)-t}\\
 & ( h=1,2,\cdots),
  \end{split}
\end{align}
where
$N_{2}^{\pm}(n,q)=|DC_{2}^{\pm}(n,q)|=A_{2}^{\pm}(n,q)B_{2}^{\pm}(n,q)$,
and $\{C_{2,j}^{\pm} (n,q)\}_{j=0}^{N_{2}^{\pm}(n,q)}$ is the weight
distribution of $C(DC_{2}^{\pm} (n,q))$ given by

\begin{align}\label{a14}
  \begin{split}
  C_{2,j}^{\pm}(n,q)=&\sum_{}^{}\binom{q^{-1}A_{2}^{\pm}(n,q)(B_{2}^{\pm}(n,q)+q^{3}-q^2-1)}{\nu_{0}}\\
                    & \times \prod_{\substack{|\tau|< 2\sqrt{q}\\ \tau\equiv -1(4)}}\;\;\prod_{K(\lambda; \beta^{-1})=\tau}\binom{q^{-1}A_{2}^{\pm}(n,q)(B_{2}^{\pm}(n,q)+q \tau -q^2 -1)}{\nu_
                    {\beta}}.
  \end{split}
\end{align}
Here the sum is over all the sets of nonnegative integers $\{\nu_
\beta\}_{\beta \in \mathbb{F}_{q}}$ satisfying  $\sum_{\beta \in
\mathbb{F}_{q}}\nu_{\beta}=j$, and $\sum_{ \beta \in \mathbb{F}_{q}
}\nu_{\beta} \beta=0$.
\end{theorem}

The following corollary is just the $n=2$ and  $n=1$ cases of  (a)
in the above. It is amusing to note that the recursive formula in
(15) and (16), obtained from the binary code $C(DC_1^-(1,q)$
associated with the double coset $DC_1^-(1,q)=P^+(2,q)$, is the same
as the one in ([5], (1), (2)), gotten from the binary code
$C(SO^+(2,q))$ associated with the special orthogonal group $SO^+
(2,q)$.

\begin{corollary}
(a) For all $q$, and $h=1,2,\cdots $,
\begin{align*}
  \begin{split}
 &MK^h=\sum_{l=0}^{h-1}(-1)^{h+l+1}\binom{h}{l}(q^2-1)^{h-l}MK^l\\
 &+q^{1-2h}(q^2-1)^{-h}\sum_{j=0}^{min \{q^2(q^2 -1)^2,h
 \}}(-1)^{h+j}C_{1,j}^{+}(2,q)\sum_{t=j}^{h}t!S(h,t)2^{h-t}\binom{q^2(q^2-1)^2-j}{q^2(q^2-1)^2-t},
  \end{split}
\end{align*}
where $\{C_{1,j}^{+}(2,q)\}_{j=0}^{q^2(q^2-1)^2}$ is the weight
distribution of $C(DC_{1}^{+}(2,q))$ given by

\begin{equation*}
C_{1,j}^+(2,q)=\sum_{}^{}\binom{q^3(q^2-1)}{\nu_0}\prod_{tr(\beta^{-1})=0}
\binom{q^2(q-1)(q+1)^2}{\nu_{\beta}}
\prod_{tr(\beta^{-1})=1}\binom{q^2(q+1)(q-1)^2}{\nu_{\beta}}.
\end{equation*}
Here the sum is over all the sets of nonnegative integers $\{\nu_
{\beta}\}_{\beta \in \fd_q }$ satisfying $\sum_{\beta \in \fd_q }
\nu_{\beta}=j$ and $\sum_{\beta \in \fd_q }\nu_ {\beta}\beta=0$. In
addition, $S(h,t)$ is the Stirling number of the second kind as
defined in (11).

(b) Let $q \geq 8 $. For $h=1,2,\cdots $,
\begin{align}\label{a15}
  \begin{split}
 MK^h=&\sum_{l=0}^{h-1}(-1)^{h+l+1}\binom{h}{l}(q-1)^{h-l}MK^l\\
      &+q\sum_{j=0}^{min \{q -1,h
      \}}(-1)^{h+j}C_{1,j}^{-}(1,q)\sum_{t=j}^{h}t!S(h,t)2^{h-t}\binom{q-1-j}{q-1-t},
  \end{split}
\end{align}
where $\{C_{1,j}^-(1,q)\}_{j=0}^{q-1}$ is the weight distribution of
$C(DC_{1}^{-}(n,q))$ given by

\begin{equation}\label{a16}
C_{1,j}^-(n,q)=\sum
\binom{1}{\nu_0}\prod_{tr(\beta^{-1})=0}\binom{2}{\nu_{\beta}}.
\end{equation}
Here the sum is over all the sets of nonnegative integers $\{\nu_0\}
\cup \{\nu_{\beta}\}_{tr(\beta^{-1})=0}$ satisfying $\nu_0+
\sum_{tr(\beta^{-1})=0}^{}\nu_ \beta=j$ and
$\sum_{tr(\beta^{-1})=0}\nu_ {\beta}\beta=0$.
\end{corollary}

\section{$O^{+}(2n,q)$}
For more details about this section, one is referred to the paper
\cite{DY}. Throughout this paper, the following notations will be
used:
\begin{itemize}
 \item [] $q = 2^r$ ($r \in \mathbb{Z}_{>0}$),\\
 \item [] $\mathbb{F}_{q}$ = the finite field with $q$ elements,\\
 \item [] $Tr A$ = the trace of $A$ for a square matrix $A$,\\
 \item [] $^tB$ = the transpose of $B$ for any matrix $B$.
\end{itemize}\

Let $\theta^+$ be the nondegenerate quadratic form on the vector
space $\fd_q^{2n \times 1}$ of all $2n \times 1$ column vectors over
$\fd_q$, given by

\begin{equation*}
\theta^{+}(\sum_{i=1}^{2n} x_i e^i) = \sum_{i=1}^{n} x_{i}x_{n+i},
\end{equation*}

where $\{e^1=^t[10\ldots0], e^2=^t[010\ldots
0],\ldots,e^{2n}=^t[0\ldots01]\}$ is the standard basis of
$\fd_q^{2n \times 1}$

The group $ O^+(2n,q)$ of all isometries of $(\fd_q^{2n \times 1},
\; \theta^+ )$ is given by:

\begin{align*}
O^{+}(2n,q)&=
\bigg\{ \left[%
\begin{smallmatrix}
  A & B \\
  C & D \\
\end{smallmatrix}%
\right]
\in GL(2n,q) \bigg| \substack{ {}^{t}AC, \;{}^{t}BD\;\; \textrm{are alternating} \\
\\
{}^{t}AD+{}^{t}CB= 1_n}\bigg\} \qquad \qquad\;\;
\end{align*}

\begin{align*}
&=
\bigg\{ \left[%
\begin{smallmatrix}
  A & B \\
  C & D \\
\end{smallmatrix}%
\right]
\in GL(2n,q) \bigg| \substack{ {}^{t}AB, \;\;{}^{t}CD \textrm{are alternating} \\
\\
A\;\;{}^{t}D+B\;\;{}^{t}C= 1_n}\bigg\},
\end{align*}
where $A, B, C, D$ are of size $n$ .

Here an $n \times n$ matrix $(a_{ij})$ is called alternating if

\begin{equation*}
\begin{cases}
 a_{ii}=0,              & \text{for $1 \leq i \leq n$},\\
 a_{ij}= -a_{ji}=a_{ji}, & \text{for $1 \leq i < j \leq n$.}
\end{cases}
\end{equation*}

$P^+=P^+(2n,q)$ is the maximal parabolic subgroup of
$O^{+}(2n,q)$ defined by:\\

\begin{align*}
P^{+}(2n,q)&=
\bigg\{ \left[%
\begin{smallmatrix}
  A & 0 \\
  0& {}^{t}A^{-1} \\
\end{smallmatrix}%
\right]
\left[%
\begin{smallmatrix}
  1_{n} & B    \\
  0     & 1_{n} \\
\end{smallmatrix}%
\right] \bigg| A \in GL(n,q),\;\; B \;\;\textrm{alternating}\bigg\}.
\end{align*}

Then, with respect to $P^+=P^+(2n,q)$, the Bruhat decomposition of
$O^+(2n,q)$ is given by:
\begin{equation}\label{a17}
O^+(2n,q)=\coprod_{r=0}^{n} P^+ \sigma_{r}^{+} P^+,
\end{equation}
where
\[
\sigma_{r}^{+}=
\begin{bmatrix}
   0   & 0         & 1_r  & 0 \\
   0   & 1_{n-r}  & 0    & 0 \\
   1_{r} & 0       & 0  &0\\
   0   & 0         & 0    & 1_{n-r}
\end{bmatrix}
\in O^+(2n,q).
\]

Put, for $0 \leq r \leq n$,
\begin{equation*}
A_{r}^{+}=\{w \in P^{+}(2n,q)| \sigma_{r}^{+}w (\sigma_{r}^{+})^{-1}
\in P^+(2n,q)\}.
\end{equation*}

Expressing $O^+(2n,q)$ as a disjoint union of right cosets of
$P^+=P^+(2n,q)$, the Bruhat decomposition in (\ref{a17}) can be
written as

\begin{equation}\label{a18}
O^+(2n,q)=\coprod_{r=0}^{n}P^{+} \sigma_{r}^{+}(A_{r}^{+}\backslash
P^{+}).
\end{equation}

The order of the general linear group $ GL(n,q)$ is given by

\begin{equation*}
g_n=\prod_{j=0}^{n-1}(q^n-q^j)=q^{\binom{n}{2}}\prod_{j=1}^{n}(q^{j}-1).
\end{equation*}

For integers $n,r$ with $0 \leq r \leq n$, the $q$-binomial
coefficients are defined as:

\begin{equation*}
\left[ \substack{n \\ r}
 \right]_q = \prod_{j=0}^{r-1} (q^{n-j} - 1)/(q^{r-j}-1).
 \end{equation*}

Then, for integers $n,r$  with $0 \leq r \leq n$, we have
\begin{equation}\label{a19}
\frac{g_n}{g_{n-r} g_r} = q^{r(n-r)}\left[ \substack{n \\ r}
 \right]_q.
\end{equation}

As it is shown in [9],
\begin{equation}\label{a20}
|A_{r}^{+}|=g_{r}g_{n-r}q^{\binom{n}{2}}q^{r(2n-3r+1)/2}.
\end{equation}
Also, it is immediate to see that
\begin{equation}\label{a21}
|P^{+}(2n,q)|=q^{\binom{n}{2}}g_n.
\end{equation}
Thus we get, from (\ref{a19})-(\ref{a21}),

\begin{equation}\label{a22}
\mid A_{r}^{+}\backslash P^{+}(2n,q) \mid = \left[ \substack{n\\
r} \right]_q q^{\binom{r}{2}},
\end{equation}
and
\begin{equation}\label{a23}
\mid P^{+}(2n,q)\sigma_{r}^{+}P^{+}(2n,q)\mid=\mid
P^{+}(2n,q)\mid^{2} \mid A_{r}^{+}
\mid^{-1}=q^{\binom{n}{2}}g_{n}\left[ \substack{n\\r} \right]_q
q^{\binom{r}{2}}.
\end{equation}

Let
\begin{equation}\label{a24}
DC_1^+(n,q)=P^+(2n,q)\sigma_{n-1}^+P^+(2n,q), \;\;\textmd{for}
\;\;n=2,4,6, \cdots,
\end{equation}
\begin{equation}\label{a25}
DC_2^+(n,q)=P^+(2n,q)\sigma_{n-2}^+ P^+(2n,q),\;\;\textmd{for}
\;\;n=2,4,6, \cdots,
\end{equation}
\begin{equation}\label{a26}
DC_1^-(n,q)=P^+(2n,q)\sigma_{n-1}^+ P^+(2n,q),\;\;\textmd{for}
\;\;n=1,3,5, \cdots,
\end{equation}
\begin{equation}\label{a27}
DC_2^-(n,q)=P^+(2n,q)\sigma_{n-2}^+P^+(2n,q),\;\;\textmd{for}
\;\;n=3,5,7, \cdots.
\end{equation}

Then, from (\ref{a23}), we have
\begin{equation}\label{a28}
N_{i}^{\pm}(n,q)=|DC_i^{\pm}(n,q)|=A_i^{\pm}(n,q)B_i^{\pm}(n,q),
\;\;\textmd{for} \;\;i=1,2
\end{equation}
(cf. (\ref{a1})-(\ref{a8})).

\textit{Unless otherwise stated, from now on, we will agree that
anything related to $ DC_{1}^{+}(n,q)$ and $DC_{1}^{-}(n,q) $ are
defined for $ n=2,4,6, \cdots$, anything related to $DC_{1}^{-}(n,q)
$ for $n=1,3,5, \cdots$, and that anything related to
$DC_{2}^{-}(n,q) $ is defined for $n=3,5,7 \cdots$.}

Also, from (\ref{a18}), (\ref{a23}), we have
\begin{align*}
  \begin{split}
|O^{+}(2n,q)|&=\sum_{r=0}^{n}|P^{+}(2n,q)|^2|A_{r}^{+}|^{-1}\\
             &=2q^{n^2-n}(q^n-1)\prod_{j=1}^{n-1}(q^{2j}-1),
  \end{split}
\end{align*}
where one can apply the following $q$-binomial theorem with $x=-1$.
\begin{equation*}
\sum_{r=0}^{n}\left[ \substack{n\\r} \right]_q(-1)^{r}
q^{\binom{r}{2}} x^r=(x; q)_{n},
\end{equation*}

with $(x;q)_n=(1-x)(1-qx)\cdots (1-q^{n-1}x)$ ($x$ an indeterminate,
$n \in \mathbb{Z}_{>0}$).

\section{Exponential sums over double cosets of $O^+(2n,2^r)$}
The following notations will be used throughout this paper.
\begin{gather*}
tr(x)=x+x^2+\cdots+x^{2^{r-1}} \text{the trace function}
~\mathbb{F}_{q}
\rightarrow \mathbb{F}_2,\\
\lambda(x) = (-1)^{tr(x)} ~\text{the canonical additive character
of} ~\fd_q.
\end{gather*}
Then any nontrivial additive character $\psi$ of $\fd_q$ is given by
$\psi(x) = \lambda(ax)$ , for a unique $a \in \fd_q^*$.

For any nontrivial additive character $\psi$ of $\fd_q$ and $a \in
\fd_q^*$, the Kloosterman sum $K_{GL(t,q)}(\psi ; a)$ for $GL(t,q)$
is defined as
\begin{equation*}
K_{GL(t,q)}(\psi ; a) = \sum_{w \in GL(t,q)} \psi(Trw + a~Trw^{-1}).
\end{equation*}
Notice that, for $t=1 $, $ K_{GL(1,q)}( \psi;a)$ denotes the
Kloosterman sum $K(\psi;a) $.

For the Kloosterman sum $K(\psi;a)$, we have the Weil bound(cf.
\cite{GJ})
\begin{equation}\label{a29}
\mid K(\psi ; a) \mid \leq 2\sqrt{q}.
\end{equation}

In \cite{D1}, it is shown that $K_{GL(t,q)}(\psi ; a)$ ~satisfies
the following recursive relation: for integers $t \geq 2$, ~$a \in
\fd_q^*$ ,
\begin{multline}\label{a30}
K_{GL(t,q)}(\psi ; a) = q^{t-1}K_{GL(t-1,q)}(\psi ; a)K(\psi
;a)\\
+ q^{2t-2}(q^{t-1}-1)K_{GL(t-2,q)}(\psi ; a),
\end{multline}
where we understand that $K_{GL(0,q)}(\psi ; a)=1$ . From
(\ref{a30}), an explicit expression of the Kloosterman sum for $GL(t,q)$ was derived  in \cite{D1}.\\

\begin{theorem}\label{C}(\cite{D1}): For integers $t \geq 1$, and $a \in \fd_q^*$, the
Kloosterman sum $K_{GL(t,q)}(\psi ; a)$ is given by
\begin{align*}
 K_{GL(t,q)}(\psi ; a)=q^{(t-2)(t+1)/2}\sum_{l=1}^{[(t+2)/2]} q^l K(\psi;a)^{t+2-2l}\sum \prod_{\nu=1}^{l-1} (q^{j_\nu -2\nu}-1),
\end{align*}
where  $K(\psi;a)$ is the Kloosterman sum and the inner sum is over
all integers $j_1,\ldots,j_{l-1}$ satisfying $2l-1 \leq j_{l-1} \leq
j_{l-2} \leq \cdots \leq j_1 \leq t+1$. Here we agree that the inner
sum is $1$ for $l=1$.
\end{theorem}

In Section 6 of \cite{DY}, it is shown that the Gauss sum for
$O^+(2n,q)$ is given by:

\begin{align}\label{a31}
\begin{split}
  \sum_{w \in O^+(2n,q)} \psi(Tr w)&=\sum_{r=0}^{n} \sum _{ w \in p^{+}\sigma_{r}^{+}P^{+}}^{} \psi(Tr w)\\
  &=\sum_{r=0}^{n}|A_{r}^{+} \backslash P^{+}| \sum_{w \in P^{+} } \psi(Tr w \sigma_{r}^{+})\\
  &=q^{\binom{n}{2}}\sum_{r=0}^{n} |A_{r}^{+} \backslash P^{+} | q^{r(n-r)}s_{r}K_{GL(n-r,q)}( \psi;1).
\end{split}
\end{align}

Here $\psi$ is any nontrivial additive character of $\fd_{q}$,
$s_{0}=1$, and, for $r \in \z_{>0}$, $s_r$ denotes the number of all
$r \times r$ nonsingular symmetric matrices over $\fd_q$, which is
given by
\begin{equation}\label{a32}
s_{r}= \begin{cases}
    q^{r(r+2)/4}\prod_{j=1}^{r/2}(q^{2j-1}-1), & \hbox{if $r$ is even,} \\
    q^{(r^2-1)/4} \prod_{j=1}^{(r+1)/2}(q^{2j-1}-1), & \hbox{if $r$ is odd,} \\
\end{cases}
\end{equation}
(cf. \cite{DY},  Proposition 4.3).

Thus we see from (\ref{a31}), (\ref{a32}), and (\ref{a22}) that, for
each $ r$ with $0 \leq r \leq n$,
\begin{equation}\label{a33}
  \begin{split}
&\sum_{w \in P^{+} \sigma_{r}^{+}P^{+}}^{}\psi(Trw)\\
&=\begin{cases}
 q^{\binom{n}{2}}q^{rn-\frac{1}{4}r^2}\left[\substack{n\\r}\right]_q \prod_{j=1}^{r/2}(q^{2j-1}-1)K_{GL(n-r,q)}(\psi;1), & \hbox{if $r$ is even,} \\
 q^{\binom{n}{2}}q^{rn-\frac{1}{4}(r+1)^2}\left[\substack{n\\r}\right]_q \prod_{j=1}^{(r+1)/2}(q^{2j-1}-1)K_{GL(n-r,q)}(\psi;1), & \hbox{if $r$ is odd.} \\
\end{cases}
\end{split}
\end{equation}

For our purposes, we need four infinite families of exponential sums
in (\ref{a33}) over $DC_{1}^{+}(n,q)$ and $DC_{2}^{+}(n,q)$ for
$n=2,4,6,\cdots$, $DC_{1}^{-}(n,q)$ for $n=1,3,5, \cdots$, and $DC_{2}^{-}(n,q)$ for $n=3,5,7,\cdots$. So we state them separately as a theorem.\\

\begin{theorem}
Let $\psi$  be any nontrivial additive character of $\fd_q$. Then,
in the notations of (\ref{a1}), (\ref{a3}), (\ref{a5}), (\ref{a7}),
we have
\begin{equation*}
\sum_{w \in DC_{1}^{\pm}(n,q)} \psi(Trw)=A_{1}^{\pm}(n,q)K(\psi;1),
\end{equation*}
\begin{align*}
  \begin{split}
  \sum_{w \in DC_{2}^{\pm}(n,q)}\psi(Trw)&=q^{-1}A_{2}^{\pm}(n,q)K_{GL(2,q)}(\psi;1)\\
                                         &=A_{2}^{\pm}(n,q)(K(\psi;1)^2+q^2-q).
  \end{split}
\end{align*}
(cf. (\ref{a33}), (\ref{a30})).\\
\end{theorem}

\begin{proposition}(\cite{D3}):\label{E}
For $n=2^s(s \in \mathbb{Z}_{\geq 0})$, and $\psi$  a nontrivial
additive character of $\fd_q$,
\[
K(\psi;a^n) = K(\psi;a).
\]
\end{proposition}
We need a result of Carlitz for the next corollary.\\

\begin{theorem}\label{F}(\cite{L2}):
For the canonical additive character $\lambda$ of $\fd_q$, and $a
\in \fd_{q} ^{*}$,
\begin{equation}\label{a34}
K_{2}(\lambda;a) = K(- \lambda;a)^{2}-q.
\end{equation}
\end{theorem}

The next corollary follows from Theorem 4, Proposition 5,
(\ref{a34}), and simple change of variables.\\

\begin{corollary}\label{G}
Let $\lambda $ be the canonical additive character of $\fd_{q}$, and
let $a \in \fd_{q}^{*}$. Then we have

\begin{align}\label{a35}
 \sum_{w \in DC_{1}^{\pm}(n,q)}\lambda(aTr w)=A_{1}^{\pm}(n,q)K(\lambda;a),
\end{align}

\begin{align}\label{a36}
  \begin{split}
 \sum_{w \in DC_{2}^{\pm}(n,q)} \lambda(aTrw)&=A_{2}^{\pm}(n,q)(K(\lambda ;a)^2+q^2-q)\\
                                      &=A_{2}^{\pm}(n,q)(K_2(\lambda;a)+q^2).
  \end{split}
\end{align}
(cf. (\ref{a1}), (\ref{a3}), (\ref{a5}), (\ref{a7})).\\
\end{corollary}

\begin{proposition}(\cite{D3})\label{H}
Let $\lambda$ be the canonical additive character of $ \fd_{q}$, $ m
\in \z_{>0}$, $\beta \in \fd_{q}$. Then

\begin{align}\label{a37}
  \begin{split}
 &\sum_{a \in \fd_{q}^{*}} \lambda(-a \beta)K_m (\lambda ;a)\\
 &=\begin{cases}
    qK_{m-1}(\lambda ;\beta^{-1})+(-1)^{m+1}, & \hbox{if $\beta \neq 0$,} \\
      (-1)^{m+1}, & \hbox{if $\beta = 0$,} \\
\end{cases}
\end{split}
\end{align}
with the convention $K_{0}(\lambda ; \beta^{-1})= \lambda(\beta
^{-1})$.
\end{proposition}

For any integer $r$ with $ 0 \leq r \leq n$, and each $\beta \in
\fd_{q}$, we let
\begin{equation*}
N_{P^{*} \sigma_{r}^{*} P^{*}}(\beta)=|\{w \in P^+ \sigma_r^+ P^+
|Tr w= \beta \}|.
\end{equation*}

Then it is easy to see that
\begin{equation}\label{a38}
qN_{P^{+}\sigma_{r}^{+}P^{+}}(\beta)=|P^{+}
\sigma_{r}^{+}P^{+}|+\sum_{a \in \fd_{q}^{*}} \lambda(-a \beta)
\sum_{w \in P^{+} \sigma_{r}^{+} P^{+} }\lambda(a Tr w).
\end{equation}
Now, from (\ref{a35})-(\ref{a38}), (\ref{a24})-(\ref{a28}), and
(\ref{a1})-(\ref{a8}), we have the following result.\\

\begin{proposition}\label{I}
\begin{align}\label{a39}
  \begin{split}
(a)\;\;&N_{DC_{1}^{\pm}(n,q)}(\beta)\\
&=q^{-1}A_{1}^{\pm}(n,q)B_{1}^{\pm}(n,q)+q^{-1}A_{1}^{\pm}(n,q)\times
\begin{cases} 1,& \beta=0,\\
q+1, & tr(\beta^{-1})=0,\\
-q+1,&tr( \beta^{-1})=1,
\end{cases}
  \end{split}
\end{align}

\begin{align}\label{a40}
  \begin{split}
(b)\;\; &N_{DC_{2}^{\pm}(n,q)}(\beta) \\
 &=q^{-1}A_{2}^{\pm}(n,q)B_{2}^{\pm}(n,q)+q^{-1}A_{2}^{\pm}(n,q)\times
 \begin{cases} qK(\lambda ; \beta^{-1})-q^{2}-1,& \beta \neq 0,\\
 q^3-q^2-1,& \beta=0.
 \end{cases}
  \end{split}
\end{align}
\end{proposition}

\begin{corollary}\label{J}
(a) For all even $n \geq 2$ and all $q$, $N_{DC_1^+(n,q)}(\beta)>
0$, for all $\beta$.

(b) For all even $ n \geq 4$ and all $ q$, or $ n=2$ and all $ q
\geq 4$, $N_{DC_{2}^{+}(n,q)}(\beta)> 0$, for all $\beta$; for $n=2$
and all $q=2$,

\begin{equation*}
N_{DC_{2}^{+}(2,2)}(\beta)
 =\begin{cases}
   0 , & \beta =1,\\
    12=|P^{+}(4,2)| ,& \beta=0.
\end{cases}
\end{equation*}
(c) For all odd $n \geq 3$ and all $q$, $N_{DC_{1}^{-}(n,q)}(\beta)>
0$, for all $\beta$; for $n=1 $ and all $q$,

\begin{equation}\label{a41}
N_{DC^{-}(1,q)}(\beta)
 =\begin{cases}
    1, & \beta =0 ,\\
     2 , & tr(\beta^{-1})=0, \\
      0, & tr(\beta^{-1})=1.
\end{cases}
\end{equation}
(d) For all odd $n \geq 3$ and all $q$, $N_{DC_2^-(n,q)}(\beta)> 0$,
for all $\beta$.
\end{corollary}

\proof (a), (c), and (d) are left to the reader.

(b) Let $n=2$. Let $\beta \neq 0$. Then, from (\ref{a40}), we have
\begin{equation}\label{a42}
N_{DC_{2}^{+}(2,q)}(0)=q^2 \{q^2-2q-1 + K(\lambda ; \beta^{-1})\},
\end{equation}

where $q^2-2q-1+ K(\lambda; \beta^{-1}) \geq q^2-2q-1-2 \sqrt{q}>0$,
for $q \geq 4$, by invoking the Weil bound in (\ref{a29}). Also,
observe from (\ref{a42}) that $N_{DC_{2}^+(2,2)}(1)=0$.

On the other hand, if $\beta =0 $, then, from (\ref{a40}), we get
\begin{equation*}
N_{DC_{2}^{+}(2,q)}(0)=q^2(2q^2-2q-1 )>0, \;\;\textmd{for all}\;\;q
\geq 2.
\end{equation*}
In addition, we note that $ N_{DC_{2}^{+}(2,2)}(0)=12$.

Assume now that $n \geq 4 $. If $\beta=0$, then, from (\ref{a40}),
we see that $N_{DC_{2}^{+}(n,q)}(0)>0$, for all $q$. Let $\beta \neq
0$. Then, again by invoking the Weil bound,
\begin{align*}
  \begin{split}
 &N_{DC_{2}^+ (n,q)}(\beta) \geq q^{-1}A_{2}^{+}(n,q)\\
 & \times \{(q^n-1)(q^{n-1}-1)q^{\frac{1}{4}(n -4)^2-1}\prod_{j=1}^{(n-2)/2}(q^{2j}-1)-(q^2 +2q^{\frac{3}{2}}+1)\}.
  \end{split}
\end{align*}

Clearly, $\prod_{j=1}^{(n-2)/2} (q^{2j}-1)>1$. So we only need to
show, for all $ q \geq 2$,
\begin{align*}
 f(q)=(q^n-1)(q^{n-1}-1)q^{\frac{1}{4}(n-4)^2-1}-(q^2+2q^{\frac{3}{2}}+1)>0.
\end{align*}
But, as $n \geq 4 $, $f(q) \geq
q^{-1}(q^4-1)(q^3-1)-(q^2+2q^{\frac{3}{2}}+1)>0$, for all $q \geq
2$.

\qquad\qquad\qquad\qquad\qquad\qquad\qquad\qquad\qquad\qquad\qquad\qquad \qquad \qquad \qquad \qquad \qquad $\square$\\

\section{Construction of codes}
Here we will construct four infinite families of binary linear codes
$C(DC_1^+(n,q))$ of length $N_1^+(n,q)$, for $n=2,4,6,\cdots$ and
all $q$, $C(DC_2^+(n,q))$ of length $N_2^+(n,q)$, for $n=2,4,6,
\cdots $ and all $q$, $C(DC_1^-(n,q))$ of length $N_1^-(n,q)$, for
$n=1,3,5,\cdots$ and all $q$, and $C(DC_2^-(n,q))$ of length
$N_2^-(n,q)$, for $n=3,5,7, \cdots$ and all $q$, respectively
associated with the double cosets $DC_1^+(n,q)$, $DC_2^+(n,q)$,
$DC_1^-(n,q)$, and $DC_2^-(n,q) $(cf. (\ref{a24})-(\ref{a27})).

Let $g_{1}, g_{2}, \cdots, g_{N_{i}^{\pm}(n,q)}$ be fixed orderings
of the elements in $DC_{i}^{\pm}(n,q)$, for $i=1,2$ by abuse of
notations. Then we put
\begin{align*}
 v_{i}^{\pm}(n,q)=(Trg_{1},Trg_{2},& \cdots,Trg_{N_{i}^{\pm}(n,q)}) \in \fd_q^{N_{i}^{\pm}(n,q)},
                             \;\;\textmd{for} \;\;i=1,2.
\end{align*}

The binary codes $C(DC_{1}^{+}(n,q))$, $C(DC_{2}^{+}(n,q))$,
$C(DC_{1}^{-}(n,q))$, and $C(DC_{2}^{-}(n,q))$ are defined as:
\begin{align}\label{a43}
 C(DC_{i}^{\pm}(n,q))=\{u \in \fd_2^{ N_{i}^{\pm}(n,q)}| u \cdot v_{i}^{\pm}(n,q)=0 \},
                                       \;\;\textmd{for}\;\; i=1,2,
\end{align}
where the dot denotes respectively the usual inner product in
$\fd_{q}^{N_{i}^{\pm}(n,q)}$, for $i=1,2$.

The following Delsarte's theorem is well-known.\\

\begin{theorem}(\cite{FN})\label{K}
Let $B$ be a linear code over $ \fd_{q}$. Then
\begin{equation*}
(B|_{\fd_{2}})^{\bot }=tr(B^{\bot}).
\end{equation*}
\end{theorem}

In view of this theorem, the respective duals  of the codes in
(\ref{a43}) are given by:
\begin{align}\label{a44}
  \begin{split}
 &C(DC_{i}^{\pm}(n,q))^{\bot}= \{c_{i}^{\pm}(a)=c_{i}^{\pm}(a;n,q)=(tr(aTrg_{1}), \cdots,tr(aTrg_{N_{i}^{\pm}(n,q)}))|a \in \fd_q \}\\
 &(i=1,2).
  \end{split}
\end{align}

Let  $\fd_2^+,\fd_q^+$ denote the additive groups of the fields
$\fd_2,\fd_q$, respectively. Then we have the following exact
sequence of groups:
\begin{equation*}
0 \rightarrow \fd_2^+ \rightarrow \fd_q^+ \rightarrow \Theta(\fd_q)
\rightarrow 0,
\end{equation*}
where the first map is the inclusion and the second one is  the
Artin-Schreier operator in characteristic two given by $x \mapsto
\Theta(x) = x^2+x$. So
\begin{equation}\label{a45}
\Theta(\fd_q) = \{\alpha^2 + \alpha \mid  \alpha \in \fd_q \},~
\textmd{and} ~~[\fd_q^+ : \Theta(\fd_q)] = 2.
\end{equation}

\begin{theorem}\label{L}(\cite{D3}):
Let $\lambda$  be the canonical additive character of $\fd_q$, and
let $\beta \in \fd_q^*$. Then
\begin{equation}\label{a46}
 (a) \sum_{\alpha \in
 \fd_q-\{0,1\}}\lambda(\frac{\beta}{\alpha^2+\alpha})=K(\lambda;\beta)-1,
 \qquad \qquad \qquad \qquad \qquad \qquad \qquad \qquad
\end{equation}
\begin{equation*}
(b)\sum_{\alpha \in
\fd_q}\lambda(\frac{\beta}{\alpha^2+\alpha+b})=-K(\lambda;\beta)-1,
\qquad \qquad \qquad \qquad\qquad \qquad \qquad \qquad
\end{equation*}
if $x^2+x+b (b\in \fd_q)$ is irreducible over $\fd_q$, or
equivalently if $b \in \fd_q\setminus\Theta(\fd_q)$ (cf.
\;(\ref{a45})).\\
\end{theorem}

\begin{theorem}\label{M}
(a) The map $\fd_q \rightarrow C(DC_{1}^+(n,q))^{\bot}(a \mapsto
c_{1}^+(a))$ is an $\fd_{2}$-linear isomorphism for $n \geq 2 $ even
and all $q $.

(b) The map $\fd_q \rightarrow C(DC_{2}^+(n,q))^{\bot}(a \mapsto
c_{2}^+(a))$ is an $\fd_{2}$-linear isomorphism for $n \geq 4 $ even
and all $q$, or $n=2$ and $q \geq 4$.

(c) The map $\fd_q \rightarrow C(DC_{1}^-(n,q))^{\bot}(a \mapsto
c_{1}^-(a))$ is an $\fd_{2}$-linear isomorphism for $n \geq 3 $ odd
and all $q $, or $n=1 $ and $q \geq 8$.

(d) The map $\fd_q \rightarrow C(DC_{2}^-(n,q))^{\bot}(a \mapsto
c_{2}^-(a))$ is an $\fd_{2}$-linear isomorphism for $n \geq 3 $ odd
and all $q$.
\end{theorem}

\proof All maps are clearly $\fd_2$-linear and surjective. Let $a$
be in the kernel of map $ \fd_q \rightarrow
C(DC_1^+(n,q))^{\bot}\;\;(a \mapsto c_1^+(a))$. Then $ tr(aTrg)=0$,
for all $g \in DC_1^+(n,q)$. Since, by Corollary 10(a), $Tr
:DC_1^+(n,q) \rightarrow \fd_q$ is surjective, $tr(a \alpha )=0$,
for all $\alpha \in \fd_q$. This implies that $a=0$, since otherwise
$ tr : \fd_q \rightarrow \fd_2$ would be the zero map. This shows
(a). All the other assertions can be handled in the same way, except
for $n=1$ and $ q \geq 8$ case of (c). Assume that we are in that
case. Then, by (\ref{a41}), $tr(a \beta)=0$, for all $\beta \in
\fd_{q}^{*}$, with $tr(\beta^{-1})=0$. Hilbert's theorem 90 says
that $tr( \gamma)=0 \Leftrightarrow \gamma=\alpha^2 + \alpha$, for
some $\alpha \in \fd_q$, and hence $\sum_{\alpha \in \fd_q -
\{0,1\}} \lambda (\frac{a}{\alpha^2+\alpha})=q-2$. If $a \neq 0$,
then, using (\ref{a46}) and the Weil bound (\ref{a29}), we would
have
\begin{equation*}
q-2=\sum_{\alpha \in \fd_q - \{0,1\}}  \lambda(\frac{a}{\alpha^2+
\alpha})=K(\lambda;a)-1 \leq 2 \sqrt{q}-1.
\end{equation*}
But this is impossible, since $x > 2 \sqrt{x} +1$, for $x \geq 8$.
\qquad \qquad  \qquad \qquad \qquad $\square$\\

Remark: One can show that the kernel of the map $\fd_q \rightarrow
C(DC_{2}^+(2,2))^{\bot}(a \mapsto c_{2}^+(a))$, and the maps $ \fd_q
\rightarrow C(DC_{1}^-(1,q))^{\bot}(a \mapsto c_{1}^-(a))$, for
$q=2,4$, are all equal to $\fd_2$.

\section{Recursive formulas for power moments of Kloosterman sums}
Here we will be able to find, via Pless power moment identity,
infinite families of recursive formulas generating power moments of
Kloosterman and 2-dimensional Kloosterman sums over all
$\fd_{q}$(with three exceptions) in terms of the frequencies of
weights in $ C(DC_{1}^+(n,q))$ or $C(DC_{1}^-(n,q))$, and
$C(DC_{2}^+(n,q))$ or $ C(DC_{2}^-(n,q))$, respectively.

\begin{theorem}\label{N}(Pless power moment identity, \cite{FN}):
Let $ B$ be an $q$-ary $[n,k]$ code, and let $B_{i}$(resp.$B_{i}
^{\bot})$ denote the number of codewords of weight $i$ in $B$(resp.
in $B^{\bot})$. Then, for $h=0,1,2, \cdots$,

\begin{align}\label{a47}
\sum_{j=0}^{n}j^{h}B_{j}=\sum_{j=0}^{min \{ n,h \}}(-1)^{j}B_{j}
^{\bot} \sum_{t=j}^{h} t! S(h,t)q^{k-t}(q-1)^{t-j}\binom{n-j}{n-t},
\end{align}
\end{theorem} where $S(h,t)$ is the Stirling number of the second
kind defined in (\ref{a11}).

\begin{lemma}
Let
\begin{equation*}
c_{i}^{\pm}(a)=(tr(a Trg_{1}), \cdots, tr(a Trg_{N_{i}^{\pm}(n,q)}))
\in C(DC_{i}^{\pm}(n,q))^{\bot},
\end{equation*} for $i=1,2$, and $a \in \fd_{q}^{*}$. Then the Hamming weights $w(c_{1}^{\pm}(a))$ and $w(c_{2}^{\pm}(a))$ are expressed
as follows:

\begin{equation}\label{a48}
(a)\;\;w(c_{1}^{\pm}(a))=
\frac{1}{2}A_{1}^{\pm}(n,q)(B_{1}^{\pm}(n,q)-K(\lambda ;
a)),\qquad\qquad \qquad\qquad
\end{equation}
\begin{equation}\label{a49}
(b)\;\;w(c_{2}^{\pm}(a))=\frac{1}{2}A_{2}^{\pm}(n,q)(B_{2}^{\pm}(n,q)-q^2+q-K(\lambda
;a)^2)\qquad\qquad
\end{equation}
\begin{equation}\label{a50}
=\frac{1}{2}A_{2}^{\pm}(n,q)(B_{2}^{\pm}(n,q)-q^2-K_2(\lambda ;
a))\qquad
\end{equation}
(cf. (\ref{a1})- (\ref{a8})).\\

\proof $w (c_{i}^{\pm}(a))=\frac{1}{2}
\sum_{j=1}^{N_{i}^{\pm}(n,q)}(1-(-1)^{tr(aTr
g_j)})=\frac{1}{2}(N_{i}^{\pm}(n,q)-\sum_{ w \in DC_{i}^{\pm}(n,q)}
\lambda(a Tr w))$, for \;$i=1,2$. our results now follow from
(\ref{a28}) and (\ref{a34})-(\ref{a36}).
\qquad \qquad \qquad \qquad \qquad $\square$\\
\end{lemma}

Let $u=(u_1, \cdots, u_{N_{N_{i}^{\pm} (n,q)}}) \in
\fd_2^{N_{i}^{\pm}(n,q)}$, for $i=1,2$, with $\nu_{\beta}$ 1's in
the coordinate places where $Tr(g_j)= \beta$, for each $\beta \in
\fd_q$. Then from the definition of the codes $C(DC_{i}^{\pm}(n,q))$
(cf. (\ref{a43})) that $u$ is a codeword with weight $j$ if and only
if $\sum_{\beta \in \fd_{q}}^{} \nu _{\beta}=j$ and $\sum_{\beta \in
\fd_{q}}^{} \nu_{\beta} \beta=0$ (an identity in $\fd_{q}$). As
there are $\prod_{\beta \in \fd_{q}}
\binom{N_{DC_{i}^{\pm}(n,q)}(\beta)}{\nu_{\beta}}$ many such
codewords with weight $j$, we obtain the following result.\\

\begin{proposition}
Let $\{C_{i,j}^{\pm}(n,q)\}_{j=0}^{N_{i}^{\pm}(n,q)}$ be the weight
distribution of $C(DC_{i}^{\pm}(n,q))$, for $i=1,2$. Then we have
\begin{equation}\label{a51}
C_{i,j}^{\pm}(n,q)=\sum \prod_{\beta \in
\fd_q}\binom{N_{DC_{i}^{\pm}(n,q)}(\beta)}{\nu_{\beta}},\textmd{
for} \;\;0 \leq j \leq N_i^{\pm}(n,q),\;\;\textmd{and}\;\;i=1,2,
\end{equation}
where the sum is over all the sets of integers $\{\nu_{\beta }\}_{
\beta \in \fd_q}(0 \leq \nu_{\beta} \leq
N_{DC_{i}^{\pm}(n,q)}(\beta))$, satisfying
\begin{equation}\label{a52}
\sum_{\beta \in \fd_{q}}^{} \nu_{\beta}=j,\;\; \textmd{and}
\sum_{\beta \in \fd_{q}}^{} \nu_{\beta} \beta=0.
\end{equation}
\end{proposition}

\begin{corollary}\label{Q}
Let $\{C_{i,j}^{\pm}(n,q)\}_{j=0}^{N_{i}^{\pm}(n,q)}$ be the weight
distribution of \;$C(DC_i^{\pm}(n,q))$, for $i=1,2$. Then we have
\begin{align*}
C_{i,j}^{\pm}(n,q)=C_{i, N_{i}^{\pm}(n,q)-j}^{\pm}(n,q),
 \;\;\textmd{for all}\;\; j,\;\; \textmd{with}\;\; 0 \leq j \leq N_{i}^{\pm}(n,q).
\end{align*}
\end{corollary}

\proof Under the replacements $\nu_{\beta}\rightarrow
N_{DC_{i}^{\pm}(n,q)}( \beta)- \nu_{\beta}$, for each  $\beta \in
\fd_q$, the first equation in (\ref{a52}) is changed to
$N_{i}^{\pm}(n,q)-j$, while the second one in there and the summands
in (\ref{a51}) are left unchanged. The second sum in (\ref{a52}) is
left unchanged, since $\sum_{\beta \in \fd_q}
N_{DC_{i}^{\pm}(n,q)}(\beta) \beta=0$, as one can see by using the
explicit expressions of $N_{DC^{\mp}(n,q)}(\beta)$ in (\ref{a39})
and (\ref{a40}).

\qquad \qquad \qquad \qquad \qquad\qquad\qquad \qquad \qquad \qquad \qquad\qquad\qquad\qquad$\square$\\

\begin{theorem}\label{R}(\cite{GJ}):
Let $q=2^r$, with $ r \geq 2$. Then the range $R$ of $K(\lambda ;a)
$, as $a$ varies over $\fd_q^{*}$, is given by:
\begin{equation*}
R=\{\tau \in \mathbb{Z} \; | \; |\tau|<2 \sqrt{q}, \; \tau \equiv -1
(mod \; 4) \}.
\end{equation*}
In addition, each value $\tau \in R $ is attained exactly $H(t^2
-q)$ times, where $H(d)$ is the Kronecker class number of $d$.
\end{theorem}

The formulas appearing in the next theorem and stated in (\ref{a10})
and (\ref{a14}) follow  by applying the formula in (\ref{a51}) to
each $C(DC_{i}^{\pm}(n,q))$, using the explicit values of
$N_{DC_{i}^{\pm}(n,q)}(\beta)$ in (\ref{a39}) and (\ref{a40}), and
taking Theorem \ref{R} into consideration.\\

\begin{theorem}\label{S}
Let $\{C_{i,j}^{\pm}(n,q)\}_{j=0}^{N_{i}^{\pm}(n,q)}$ be the weight
distribution of  $C(DC_{i}^{\pm}(n,q))$, for $i=1,2$, and assume
that $q \geq 4$, for $C(DC_2^{\pm}(n,q))$. Then we have

(a) For $j=0,\cdots,N_{1}^{\pm}(n,q)$,
\begin{align*}\label{}
  \begin{split}
 &C_{1,j}^{\pm}(n,q)=\sum \binom{q^{-1}A_{1}^{\pm}(n,q)(B_{1}^{\pm}(n,q)+1)}{\nu_0}\\
 &\times \prod_{tr(\beta^{-1})=0}\binom{q^{-1}A_{1}^{\pm}(n,q)(B_{1}^{\pm}(n,q)+q+1)}{\nu_{\beta }}\prod_{tr(\beta^{-1})=1} \binom{q^{-1}A_{1}^{\pm}(n,q)(B_{1}^{\pm}(n,q)-q+1)}{\nu_{\beta}},
  \end{split}
\end{align*}
where the sum is over all the sets of nonnegative integers $\{\nu_
{\beta}\}_{\beta \in \fd_q}$ satisfying $\sum_{\beta \in \fd_q}\nu_
\beta=j$ and $\sum_{\beta \in \fd_q}\nu_{\beta}\beta=0$.

(b) For $j=0,\cdots,N_2^{\pm}(n,q)$,
\begin{align*}
  \begin{split}
 &C_{2,j}^{\pm}(n,q)= \sum \binom{q^{-1}A_{2}^{\pm}(n,q)(B_{2}^{\pm}(n,q)+q^3-q^2-1)}{\nu_0}\\
 & \times \prod_{\substack{|\tau |< 2\sqrt{q}\\ \tau \equiv -1(4)}}\;\;\prod_{K(\lambda ; \beta^{-1})= \tau}
\binom{q^{-1}A_{2}^{\pm}(n,q)(B_{2}^{\pm}(n,q)+q \tau -q^2 -1)}
{\nu_ {\beta}},
  \end{split}
\end{align*}
where the sum is over all the sets of nonnegative integers $\{\nu_
{\beta}\}_{\beta \in \fd_q}$ satisfying $\sum_{\beta \in \fd_{q}}
\nu_{\beta}=j$, and $\sum_{\beta \in \fd_q} \nu_{\beta} \beta=0$.
\end{theorem}
From now on, we will assume that, for $C(DC_1^+(n,q))^{\bot}$, $n
\geq 2$ even and all $q$; for $C(DC_2^+(n,q))^{\bot}$, $n \geq 2$
even and $q \geq 4 $; for $C(DC_1^-(n,q))^{\bot}$, either $n \geq 3
$ odd and all $q$ , or $n=1$ and $q \geq 8$ ; for $
C(DC_2^-(n,q))^{\bot}$, $ n \geq 3$ odd and $q \geq 4$. Under these
assumptions, each codeword in $C(DC_i^{\pm}(n,q))^{\bot}$ can be
written as $c_i^{\pm}(a)$, for $i=1,2$, and a unique $a \in
\fd_q$(cf. Theorem 13, (\ref{a44})).

Now, we apply the Pless power moment identity in (\ref{a47}) to $
C(DC_{i}^{\pm}(n,q))^{\bot}$, for those values of $n$ and $q$, in
order to get the results in Theorem \ref{A} (cf. (\ref{a9}),
(\ref{a12}),(\ref{a13})) about recursive formulas.

The left hand side of that identity in (\ref{a47}) is equal to
\begin{equation*}
\sum_{a \in \fd_{q}^{*}}w(c_{i}^{\pm}(a))^h,
\end{equation*}
with $w(c_{i}^{\pm}(a))$ given by (\ref{a48})-(\ref{a50}). We have

\begin{equation*}
\sum_{a \in
\fd_{q}^{*}}w(c_{1}^{\pm}(a))^{h}=\frac{1}{2^h}A_{1}^{\pm}(n,q)^h
\sum _{a \in \fd_{q}^{*}}(B_{1}^{\pm}(n,q)-K(\lambda ; a))^h
\end{equation*}
\begin{equation}\label{a53}
=\frac{1}{2^h}A_{1}^{\pm}(n,q)^h \sum_{l=0}^{h}(-1)^l
\binom{h}{l}B_{1}^{\pm}(n,q)^{h-l} M K^l.
\end{equation}
Similarly, we have
\begin{align}\label{a54}
\sum_{a \in
\fd_q^*}w(c_{2}^{\pm}(a))^h=\frac{1}{2^h}A_{2}^{\pm}(n,q)^h
\sum_{l=0}^{h}(-1)^l \binom{h}{l}(B_{2}^{\pm}(n,q)-q^2+q)^{h-l} M
K^{2l}
\end{align}
\begin{align}\label{a55}
=\frac{1}{2^h}A_{2}^{\pm}(n,q)^h  \sum_{l=0}^{h}(-1)^l
\binom{h}{l}(B_{2}^{\pm}(n,q) -q^2)^{h-l} M K_{2}^{l}.
 \end{align}

Note here that, in view of (\ref{a34}), obtaining power moments of
2-dimensional Kloosterman sums is equivalent to getting even power
moments of Kloosterman sums. Also, one has to separate the term
corresponding to $l=h$ in (\ref{a53})-(\ref{a55}), and notes $
dim_{\fd_2}C(DC_{i}^{\pm}(n,q))^{\bot}=r$.

\begin{flushleft}
{\textbf{Acknowledgments}}
\end{flushleft}
``This work was supported by the Korea Research Foundation Grant
funded by the Korean Government(MOEHRD, Basic Research Promotion
Fund) (KRF-2008-C00007-40003038)."

\end{document}